\documentclass[]{interact}
\usepackage[numbers,sort&compress]{natbib}
\bibpunct[, ]{[}{]}{,}{n}{,}{,}
\renewcommand\bibfont{\fontsize{10}{12}\selectfont}
\makeatletter
\def\NAT@def@citea{\def\@citea{\NAT@separator}}
\makeatother
\usepackage{tabularx}
\usepackage{booktabs}
\usepackage{color}
\usepackage{caption}
\usepackage{adjustbox}
\usepackage{multirow}
\usepackage{bbm}
\usepackage{longtable} 
\usepackage{pdflscape} 
\usepackage{geometry}
\usepackage{fix-cm}
\usepackage{multicol}
\usepackage{adjustbox}
\usepackage{xr-hyper}
\usepackage[pdftex,colorlinks=true,allcolors=blue]{hyperref}

\usepackage[section]{placeins}
\usepackage{caption}
\usepackage{graphicx} 
\usepackage[font=small,skip=10pt]{caption}
\usepackage[font=footnotesize,skip=2pt]{subcaption}
\usepackage{subcaption}
\allowdisplaybreaks
\usepackage{booktabs}
\usepackage{array}
\usepackage{multirow}
\usepackage{makecell}
\usepackage{floatpag}
\usepackage{amsmath,amssymb,amsthm}
\usepackage{mathtools}
\usepackage{physics}
\usepackage{enumitem}
\usepackage{textcomp}
\usepackage{gensymb}
\usepackage{arydshln}
\usepackage{url}
\usepackage{fullpage}
\usepackage{bm}
\usepackage{algorithm}
\usepackage[noend]{algpseudocode} 

\usepackage{bookmark}
\usepackage{lmodern}

\newtheorem{definition}{Definition}

\newtheorem{theorem}{Theorem}

\newtheorem{lemma}{Lemma}

\newtheorem{proposition}{Proposition}

\usepackage{comment}

\usepackage{rotating}
\usepackage{pdflscape}
\usepackage{cleveref}

\renewcommand{\Cref}[1]{\ref{#1}}
\usepackage{float}
\newcolumntype{Y}{>{\centering\arraybackslash}X}

\usepackage{lineno}

\newcommand{\st}{\mbox{s.t.}} 

\newcounter{savealgorithm}

\newcounter{procedure}
\makeatletter

\floatpagestyle{plain}
\newcommand{\StatexIndent}[1][3]{%
  \setlength\@tempdima{\algorithmicindent}%
  \Statex\hskip\dimexpr#1\@tempdima\relax}
\algdef{S}[WHILE]{WhileNoDo}[1]{\algorithmicwhile\ #1}%
\makeatother

\hypersetup{colorlinks=true,allcolors=blue}

\newcommand{\Ltwo}{L^2(0,L)}
\newcommand{\Lomega}{L^2(\Omega)}
\newcommand{\Vad}{{\cal{V}}_{ad}}
\newcommand{\Homega}{{\hat{H}}^{2,1}(\mathrm{\Omega})}

\usepackage{lineno}
\begin{document}

\articletype{ARTICLE TEMPLATE}

\title{On the optimal control of initial velocity in a hyperbolic beam equation by the variational method}

\author{
\name{Ye\c{s}im Akbulut\textsuperscript{a}
and Bismark Singh\textsuperscript{b}}
\thanks{\textsuperscript{a}Department of Mathematics, Atat\"{u}rk University, 25240 Erzurum, Turkey. Email: \url{ysarac@atauni.edu.tr}
\newline 
\textsuperscript{b}School of Mathematical Sciences, University of Southampton, Southampton SO17 1BJ, UK. Email: \url{b.singh@southampton.ac.uk}}
}
\maketitle

\begin{abstract}
We study the problem of controlling the initial condition of a vibrating beam. The optimal control problem seeks to determine solutions of initial velocity that assure the approach of the state of the beam to a given target function in the $L^2-$norm. We prove both the existence and uniqueness of the optimal solution. Employing identities based on the adjoint and difference problems, we determine the Fr\'echet derivative of the cost functional. We further derive the necessary optimality conditions of this control problem. Finally, we provide a sketch of a gradient-based algorithm, that rests on the explicit formula of the gradient of the cost functional, to obtain numerical solutions. 
%
%
%
\end{abstract}

\begin{keywords}
Beam equations, Adjoint approach, Fr\'echet derivative, Optimal control.
\end{keywords}

\section{Preliminaries} \label{sec:notation}
We begin with a few preliminaries for  this work. Given a function $f(x,t)$, we employ the following notation: 
\begin{itemize}
 \item $f(a,b)$ denotes the value of the function $f$ at $x=a, t=b$.
 \item  $f_{xx}$ denotes the second partial derivative of $f$ with respect to its argument $x$, and similarly for $t$.
 \item $f_{x}(a,b)$ denotes the value of the first partial derivative of $f$ with respect to $x$ at the point $x=a, t=b$, and similarly for $t$.
 \item $f_{xx}(a,b)$ denotes the value of the second partial derivative of $f$ with respect to $x$ at the point $x=a, t=b$, and similarly for $t$.
 \item  For a function $g(x)$, we employ the following notation: $g'$ and $g''$ denote the first and second ordinary derivatives of $g$.
\end{itemize}
 We let $\Omega$ denote the set $(0,L)\times (0,T)$. We employ the below listed spaces and their corresponding norms.  

\begin{definition} \label{def:ltwo}
\cite{Adams1975}
 The spaces $\Ltwo$ and $\Lomega$ are Hilbert spaces consisting of all Lebesgue-measurable functions $f$ and $g$ on $(0,L)$ and $\Omega$ that have finite norms given by 
 \begin{eqnarray*}
  & \norm{f}_{\Ltwo} = & \Big(\int_0^L f(x)^2 dx\Big)^{1/2}, \\
  & \norm{g}_{L^2(\Omega)} = & \Big(\int_0^T \int_0^L g(x,t)^2 dx dt \Big)^{1/2},
 \end{eqnarray*}
 respectively. 
 Further, the inner product of two functions ${<f_1,f_2>}_{\Ltwo}$ and ${<g_1,g_2>}_{\Lomega}$ is defined as $\displaystyle  \int_0^L f_1(x) f_2(x) dx$ and $\displaystyle \int_0^T \int_0^L g_1(x,t) g_2(x,t) dx dt$ in these two spaces, respectively. \qed
\end{definition}

\begin{definition} \label{def:Linfinity}
\cite{Adams1975}
 The space $L^\infty(0,L)$  consists of equivalence classes of functions $f:(0,L)\rightarrow \mathbb{C}$ such that there is a finite constant $M$ with $\mid f(x) \mid \le M$ for almost every $x \in (0,L)$. The norm in this space is defined by
$$\norm{f}_{L^\infty(0,L)} =  \operatorname{ess}\sup\limits_{} \{|f(x)|, x\in (0,L)\}.$$
\end{definition}

\begin{definition} \label{def:Hsquare}
\cite{Adams1975}
 The space $H^2(0,L)$ is a Hilbert space consisting of all functions $f$ in $\Ltwo$ that have weak partial derivatives $f_x, f_{xx} \in \Ltwo$.  The scalar product of functions $f$ and $g$, and the norm in this space are defined by
 \begin{eqnarray*}
 & <f,g>_{H^2(0,L)} = &  {\int_0^L {\bigg[fg+f_xg_x+f_{xx}g_{xx}\bigg]dx}} \\
 & \norm{f}_{H^2(0,L)} = &\Big({<f,f>}_{H^2(0,L)}\Big)^{\frac{1}{2}},
\end{eqnarray*}
respectively. \qed
\end{definition}

\begin{definition} \label{def:htwo}
\cite{Ladyzhenskaya2013}
The space $H^{2,1}(\Omega)$ is a Hilbert space consisting of all functions $f$ in $\Lomega$ that have weak partial derivatives $f_t, f_x,f_{xx} \in \Lomega$. Then, a norm in this space is defined as:
\begin{equation*}
\norm{f}_{H^{2,1}(\Omega)}= {\Big(\int_0^T \int_0^L {\Big[f^2+f^2_t+f^2_x+f^2_{xx}\Big]dxdt}\Big)}^{\frac{1}{2}}. \hfill 
\end{equation*} \qed
\end{definition}

\begin{definition}\label{def:htwohat}
Consider Definition~\ref{def:htwo}. The space ${\hat{H}}^{2,1}(\Omega)$ is a  subset of the space $H^{2,1}(\Omega)$ consisting of all functions $f$ such that $f(0,t)=f(L,t)=0, \forall t\in [0,T]$. Then, a norm in this space is defined as in Definition~\ref{def:htwo}. \qed
\end{definition}

\section{Introduction} \label{sec:introduction}
Consider the following hyperbolic problem governing the vibrations of a beam:
\begin{subequations}\label{eqn:euler_bernoulli}
\begin{align}
& u_{tt}+\big(k(x)u_{xx}\big)_{xx}=F(x,t),  & \forall  (x,t) \in \Omega 
\label{eqn:euler_bernoulli_1} \\
& \, \,  \text{with initial conditions:}  \nonumber\\
& u(x,0)=w(x), \, u_{t}(x,0)=v(x), &\forall x\in (0,L) \label{eqn:euler_bernoulli_2} \\
& \, \,  \text{and boundary conditions:}  \nonumber\\
& u(0,t)=u(L,t)= u_{xx}(0,t)=u_{xx}(L,t)=0,                 & \forall t \in (0,T). \label{eqn:euler_bernoulli_3}
\end {align}
\end{subequations}

Model~\eqref{eqn:euler_bernoulli} is an initial-boundary value problem for the linear Euler-Bernoulli beam equation that is well-studied in the  engineering literature, see, e.g., ~\cite{Inman1994,Clough1993}. Here, the unknown functions $u$ and $v$ denote the deflection of the beam and its initial velocity, respectively. Specifically, 
\begin{itemize}
 \item $u(x,t)$ denotes the vertical displacement of the beam;
 \item $u_t(x,t)$ represents the velocity of the beam; and
 \item $u_{xx}(x,t)$ captures the curvature of the beam at a position $x \in (0,L)$ and time $t \in (0,T)$.
\end{itemize}
The model paramters are:
\begin{itemize}
 \item  $T>0$ and $ L> 0$, denoting the time horizon and beam length,  respectively;
 \item $k(x)=E \times I(x)>0$, representing the flexural rigidity with a Young's modulus $E$ and 
moment of inertia $I(x)$; and,
\item $F(x,t)$ and $w(x)$ denoting the external load and the initial displacement of the beam, respectively. 
\end{itemize}
Given a function $v(x)$, and the known functions $F(x,t), w(x)$ and $k(x)$, model~\eqref{eqn:euler_bernoulli} seeks a solution for the function $u(x,t)$. This problem is known as the \textit{direct problem} or the \textit{forward problem}.  

Unlike the direct problem, we consider a setting where the function $v(x)$ is unknown and needs to be identified alongside $u(x,t
;v)$; here, $u(x,t
;v)$ denotes the solution of the direct problem~\eqref{eqn:euler_bernoulli} corresponding to the given $v$.
Specifically, our goal is to find a suitable-sized optimal control velocity, $v(x)$, such that the beam's displacement, $u(x,t;v)$ closely approximates a given target function $y(x,t) \in \Lomega$. By closely, we mean that the quantity $\norm{u\left(x,t;v\right)-y\left(x,t\right)}_{\Lomega}$ is small; for a definition of a norm in $\Lomega$, see Definition~\ref{def:ltwo}. To do so, we restrict ourselves to a set of \textit{admissible control} functions belonging to the space $\Ltwo$ (see, Definition~\ref{def:ltwo}) as follows:
 \begin{equation}\label{eqn:vad}
{\cal{V}}_{ad}=\{v\in \Ltwo: \norm{v}_{\Ltwo} \le v_c \},
\end{equation}
where $v_c$ is a given finite bound. 
Thus, our optimal control problem is formulated as follows:
\begin{subequations}\label{eqn:definition_j}
\begin{align}
 J^*_{\alpha} & =   \inf_{v \in \Vad}  J_{\alpha}(v) \label{eqn:definition_j_obj} \\
 &  \st \,  J_{\alpha}(v)=  \displaystyle  \int^T_0{\int^L_0{{\big[u\left(x,t;v\right)-y\left(x,t\right)\big]}^2dx}dt}+\alpha \int^L_0{v^2\left(x\right)dx}.\label{eqn:definition_j_constraint} 
\end{align}
\end{subequations}

In model~\eqref{eqn:definition_j}, the function $y(x,t)$ is the given target function and $\alpha >0$ is a regularization parameter. The objective function~\eqref{eqn:definition_j_obj} seeks  to identify an {optimal} control function, $v^*$, from the set of admissible controls, $\Vad$, that minimizes the cost functional, $J_{\alpha}(v)$, defined by equation~\eqref{eqn:definition_j_constraint}. The first term in~\eqref{eqn:definition_j_constraint} represents the approximation error, while the second term serves as a regularization or penalty term for stabilization.  With this background, the following model presents the the central problem of our work.
\begin{subequations}\label{eqn:main_problem}
\begin{align}
  J^*_{\alpha} & = \inf_{v \in \Vad} J_{\alpha}(v) & \label{eqn:main_problem_obj} \\
   &  \st \,  J_{\alpha}(v)=  \displaystyle  \int^T_0{\int^L_0{{\big[u\left(x,t;v\right)-y\left(x,t\right)\big]}^2dx}dt}+\alpha \int^L_0{v^2\left(x\right)dx}   \label{eqn:main_problem_1} & \\
   & u_{tt}+\big(k(x)u_{xx}\big)_{xx}=  F(x,t),  & \forall  (x,t) \in \Omega  \label{eqn:main_problem_2} \\
 & \quad \text{with initial conditions:} &   \nonumber \\
&  u(x,0)  =  w(x), \, u_{t}(x,0)=v(x), &  \forall x\in (0,L) \label{eqn:main_problem_3} \\
 & \quad \text{with boundary conditions:} &    \nonumber\\
& u(0,t)=u(L,t)=  0,  u_{xx}(0,t)=u_{xx}(L,t)=0, & \forall t \in (0,T). \label{eqn:main_problem_4} 
\end{align}
\end{subequations}

By solving model~\eqref{eqn:main_problem}, we obtain an optimal control function $v^*$ that correspondingly determines $u(x,t;v)$ via model~\eqref{eqn:main_problem_2}-\eqref{eqn:main_problem_4} such that it is close enough to the given target function $y(x,t)$.  

The Euler-Bernoulli beam equation has been extensively studied in engineering and the applied sciences. It is frequently used to formulate problems describing  the vibration of beams that find application in  structural mechanics, aerospace engineering, and mechanical engineering. For foundational studies on beam vibration problems, see,  e.g.,~\cite{Meirovitch2010,Popov1978}. Recent studies on the Euler-Bernoulli beam equations focus on inverse problems with uncertain or missing data. Such problems are often posed as optimal boundary control problems, see, e.g.,~\cite{Conrad1998, Guo2002,Hasanov2019,Karagiannis2018}. Unlike these works, we consider the initial velocity $v(x)$ as a control variable. 

Although the optimal control of initial conditions in a hyperbolic system is studied~\cite{Sarac2012}, there are relatively few works on controlling the initial conditions in beam equations.  Here, two works stand out as notable exceptions relevant to this work. First, 
First, Engin et al.~\cite{Engin2022} study a similar problem as model~\eqref{eqn:main_problem} with the difference that their initial  velocity, $v$, is controlled via the beam's final displacement, $u(x,T)$, and final velocity, $u_t(x,T)$; in contrast, our  problem controls the initial velocity by approximating the beam's vertical displacement, $u(x,t)$, to a target function. Second,  Hwang~\cite{Hwang2023} studies a minimax optimal control problem  of a nonlinear beam equation with distributed controls and noise in the initial velocity and derives the necessary optimality conditions via the Pontryagin's maximum principle. 

The rest of this work is organized as follows. In Section~\ref{sec:existence}, we show both the existence and the uniqueness of a weak solution $u$, for a given function $v$, of model~\eqref{eqn:euler_bernoulli} as well as both the existence and the uniqueness of the optimal solution, $v^*$, of model~\eqref{eqn:main_problem}. Building on this, Section~\ref{sec:frechet} provides the necessary optimality conditions using an adjoint-based approach.  In Section~\ref{sec:lipschitz}, we show the Lipschitz continuity of the cost functional's gradient that can be used for gradient-based numerical solutions of model~\eqref{eqn:main_problem}. We conclude with a summary of our work, as well as directions for future work, in Section~\ref{sec:conclusion}.

\section{Existence and uniqueness of an optimal solution for model~\eqref{eqn:main_problem}} \label{sec:existence}

We begin by analyzing a solution, $u$, of model~\eqref{eqn:euler_bernoulli} in the weak sense as proposed by~\cite{Ladyzhenskaya2013}. 

\begin{definition} \label{def:weak_solution}
 Consider Definition~\ref{def:htwohat} and  model~\eqref{eqn:euler_bernoulli}. Then, given a function $v$ from the admissible control set $\Vad$ defined by equation~\eqref{eqn:vad}, a \textit{weak solution} of model~\eqref{eqn:euler_bernoulli} is a function $u$ that satisfies: (i) $u\in \Homega$, and (ii) $u(x,0) = w (x), \forall x \in (0,L)$ if
 \[
 \int^T_0{\int^L_0{\big(-u_t{\eta }_t+k(x)u_{xx}{\eta }_{xx}\big)dxdt}}=\int^T_0{\int^L_0{\eta F dxdt}} +\int^L_0{v(x)\eta(x,0)dx}
\]
holds for all functions $\eta \in \Homega$ with $\eta \big(x,T\big)=0, \forall x \in (0,L)$. \qed
\end{definition}

Unless otherwise stated, in what follows, we always a refer to a solution of a system of differential equations as being in the weak sense. The following lemma establishes a sufficient condition for both the existence and uniqueness of a weak solution for model~\eqref{eqn:euler_bernoulli}.

\begin{lemma} \label{lem:weak}  Consider Definition~\ref{def:ltwo}-Definition~\ref{def:htwohat}.
\begin{enumerate}[label={L.(\roman*)}] 
\item  \label{lem:weak_a} Existence and uniqueness. 

\noindent There exists a weak solution for model~\eqref{eqn:euler_bernoulli} in the sense of Definition~\ref{def:weak_solution}  if:  (a)  the function $k\big(x\big)\in L_{\infty }(0,L)$ is positive and finite-valued, (b) $w(x) \in H^2(0,L)$, and (c) $F(x,t)\in L^2(\Omega)$. Further, this weak solution is unique.
 \item \label{lem:weak_b} Energy estimate. 
 
 \noindent Let the above conditions hold.  Then, we have 
\[ \norm{u}_{H^{2,1}(\Omega)} \le C_0 \bigg(\norm{w''}_{\Ltwo}+\norm{v}_{L^{2}(0,L)}+\norm{F}_{L^{2}(\Omega)}\bigg).\]
\end{enumerate}
\end{lemma}
\begin{proof} \,
 \begin{enumerate}[label=(\roman*)]
  \item The proof of existence follows directly using the Galerkin method as in~\cite{Ladyzhenskaya2013}. For uniqueness, assume that $u_1$ and $u_2$ are two solutions of model~\eqref{eqn:euler_bernoulli}. Then the difference function $r(x,t)=u_1(x,t)-u_2(x,t)$ satisfies the following model:
  \begin{subequations}\label{eqn:function_r}
\begin{align}
& r_{tt}+\big(k(x)r_{xx}\big)_{xx}=0,  & \forall  (x,t) \in \Omega 
\label{eqn:function_r_1} \\
& \, \,  \text{with initial conditions:}  \nonumber\\
& r(x,0)=0, \, r_{t}(x,0)=0, &\forall x\in (0,L) \label{eqn:function_r_2} \\
& \, \,  \text{and boundary conditions:}  \nonumber\\
& r(0,t)=r(L,t)= r_{xx}(0,t)=r_{xx}(L,t)=0,                 & \forall t \in (0,T). \label{eqn:function_r_3}
\end {align}
\end{subequations}
The energy estimate for model~\eqref{eqn:function_r} (obtained in a similar way as in the proof of Lemma~\ref{lem:weak_b} below) is:
\[ \norm{r}_{H^{2,1}(\Omega)} \le 0.\]
Hence, $r(x,t)$ is identically zero over $\Homega$; i.e., $u_1=u_2$ which shows the uniqueness of the solution. 
  \item The proof is similar to that given in~\cite{Hasan2009} with the exception of the boundary condition given by~\eqref{eqn:euler_bernoulli_3}. Consider the {energy} identity:
  \begin{equation} \label{eqn:lemma_weak_1}
   \big(k(x)u_{xx}\big)_{xx} u_t = \bigg[(k(x)u_{xx})_{x} u_t-k(x)u_{xx} u_{tx}\bigg]_x+\frac{k(x)}{2}\frac{d}{dt}(u_{xx}^2).
  \end{equation}
We then have the following. 

\begin{subequations} \label{eqn:lemma_weak_2}
\begin{align}
        F u_t &= \frac{1}{2} \frac{d}{dt} (u_t^2) + \frac{k(x)}{2} \frac{d}{dt} (u_{xx}^2) 
        + \bigg[(k(x) u_{xx})_x u_t - k(x) u_{xx} u_{tx} \bigg]_x, \label{eqn:lemma_weak_2_1} \\
        \implies \quad \int_0^L F u_t \, dx &= \frac{1}{2} \frac{d}{dt} \int_0^L \big[ u_t^2 + k(x) u_{xx}^2 \big] dx, \label{eqn:lemma_weak_2_2} \\
        \implies \quad \int_0^t \int_0^L F(x,\tau) u_\tau(x,\tau) \, dx \, d\tau &= 
        \frac{1}{2} \int_0^L \big[ u_t(x,t)^2 + k(x) u_{xx}(x,t)^2 \big] dx \nonumber \\
        &\quad - \frac{1}{2} \int_0^L \big[ v(x)^2 + k(x) w''(x)^2 \big] dx. \label{eqn:lemma_weak_2_3}
\end{align}
\end{subequations}

Here, equation~\eqref{eqn:lemma_weak_2_1} follows from equation~\eqref{eqn:lemma_weak_1} and by multiplying equation~\eqref{eqn:euler_bernoulli_1} by $u_t$.  Then, equation~\eqref{eqn:lemma_weak_2_2} follows by integrating equation~\eqref{eqn:lemma_weak_2_1}  on $[0,L]$; here, the integral of the term inside the square parentheses in equation~\eqref{eqn:lemma_weak_2_1} is zero from the boundary conditions~\eqref{eqn:euler_bernoulli_3}. Next, from the initial condition \eqref{eqn:euler_bernoulli_2} we have $u_{xx}(x,0)=w''(x)$. Then, equation~\eqref{eqn:lemma_weak_2_3} holds by further integrating equation~\eqref{eqn:lemma_weak_2_2} on $[0,t]$ and using the initial condition \eqref{eqn:euler_bernoulli_2}.

Equation~\eqref{eqn:lemma_weak_2_3} is the so-called energy identity that is also well-studied, see, e.g.,~\cite{Hasan2009,Hasan2024}. The weak solution of  model \eqref{eqn:euler_bernoulli} satisfies this energy identity. Further, an energy norm is defined as follows.
\begin{equation} \label{eqn:energy_norm}
 \norm{u}_E=\Big(\int_0^T\int_0^L \bigg[( u_{t}(x,t)^2+k(x)( u_{xx}(x,t)^2\bigg]dxdt\Big)^\frac{1}{2}.
\end{equation}

Then, from~\eqref{eqn:lemma_weak_2_3} and using the processes done in~\cite{Ladyzhenskaya2013} , we obtain the following estimate.
\begin{equation} \label{eqn:energy_estimate}
 \norm{u}_{E} \le C_{00} \bigg(\norm{w''}_{\Ltwo}+\norm{v}_{L^{2}(0,L)}+\norm{F}_{L^{2}(\Omega)}\bigg)
\end{equation}
\noindent Since the norms $\norm{\cdot}_{H^{2,1}(\Omega)}$ and $\norm{\cdot}_{E}$ are equivalent, the proof is complete.
\end{enumerate}
\end{proof}

\begin{theorem} \label{cor:existence}
 An optimal control, $v^*$, for model~\eqref{eqn:main_problem} exists and is unique.
\end{theorem}
\begin{proof}
 For existence, we need to show: (i) the existence of a solution for equations~\eqref{eqn:main_problem_2}-\eqref{eqn:main_problem_4}, and (ii) the cost functional~\eqref{eqn:main_problem_obj} makes sense at any such solution.
 
 Lemma~\ref{lem:weak_a} shows the existence of a unique weak solution for model~\eqref{eqn:euler_bernoulli}; hence, part~(i) follows directly since the constraints of model~\eqref{eqn:main_problem} are given by model~\eqref{eqn:euler_bernoulli}.  For part~(ii), consider a sequence $\left\{v_n\right\}\subset V_{ad}$  that converges to a  function $v\in V_{ad}$. Such a sequence always exists since the set $\Vad$ is closed and bounded. Then, from Lemma~\ref{lem:weak_b} the sequence of traces $\left\{u\left(x,t;v_n\right)\right\}$ also
converges to a function $\left\{u\left(x,t;v\right)\right\}$ in the $L^2(\Omega)$ space. From the definition of the cost functional in equation~\eqref{eqn:main_problem_obj},  a corresponding sequence of cost functionals for this trace also converges; thus, the cost functional is continuous.  Then, a minimizer for  model~\eqref{eqn:main_problem} exists by the generalized Weierstrass theorem, see, e.g.,~\cite{Zeidler1985}. 

Since the cost functional is the sum of a convex and a strictly convex function, the minimizing function is unique. This proves the result.
\end{proof}

\section{A necessary condition to verify optimality of model~\eqref{eqn:main_problem}} \label{sec:frechet}
In Section~\ref{sec:existence}, we show that a unique solution, $v^*$, of model~\eqref{eqn:main_problem} exists. Next, we derive an optimality condition that an optimal control, $v^*$,  of model~\eqref{eqn:main_problem} must satisfy.  We do so in Theorem~\ref{thm:necessarycondition} which is the central result of this section. To prove this result, we employ the Fr\'echet gradient of the cost functional that we compute using an adjoint-based approach. This in turn rests on three results that we begin with.

First, consider an increment $\delta v$ to $v$ such that $v+\delta v\in V_{ad}$, and the following problem:
\begin{subequations}\label{eqn:difference}
\begin{align}
&\delta u_{tt}+(k(x)\delta u_{xx})_{xx}=0,  & \forall (x,t) \in \Omega \label{eqn:difference_1}\\
& \, \, \text{with initial conditions:}  \nonumber\\
&\delta u(x,0)=0, ~\delta u_{t}(x,0)=\delta v(x), & \forall x \in (0,L) \label{eqn:difference_2}\\
& \, \,  \text{with boundary conditions:} \nonumber \\
&\delta u(0,t)=\delta u(L,t)=0, \delta u_{xx}(0,t)=\delta u_{xx}(L,t)=0,                 & \forall t \in (0,T). \label{eqn:difference_3}
\end {align}
\end{subequations}

Model~\eqref{eqn:difference} is known as the \textit{difference problem}. It has a solution in the space $\Homega$ given by  $ \delta u (x,t;v) = u(x,t;v+\delta v)-u(x,t;v)$.   Since model~\eqref{eqn:difference} is of the same type as model~\eqref{eqn:euler_bernoulli}, it follows from  Lemma~\ref{lem:weak_a} that this solution is also unique. Further, we have the following proposition.  

\begin{proposition}\label{cor:frechet}
 Let $\delta u(x,t;v) = u(x,t;v+\delta v)-u(x,t;v)$ denote the weak solution of model~\eqref{eqn:difference}. Then the following estimate is valid
\begin{equation}\label{e32}
\norm{\delta u}_{\Lomega} \le C_1 \norm{\delta v}_{\Ltwo},
\end{equation}
where the parameter $C_1$ is independent of $\delta v$.
\end{proposition}
\begin{proof}
 From Lemma~\ref{lem:weak_b},  we have $\norm{\delta u}_{H^{2,1}(\Omega)} \le C_0 \norm{\delta v}_{\Ltwo}$. Since $\norm{\delta u}_{\Lomega} \le \norm{\delta u}_{H^{2,1}(\Omega)}$, the result follows.
\end{proof}

\noindent Second, consider the Lagrangian   functional, $L(u,v,\psi )$,  for model~\eqref{eqn:main_problem} defined as follows. 
\begin{equation} \label{9} 
L(u,v,\psi )=J_{\alpha }(v)+{\big \langle \psi ,u_{tt}+{(k(x)u_{xx})}_{xx}-F(x,t) \big\rangle}_{L^2(\mathrm{\Omega })}, 
\end{equation} 
where the function $\psi = \psi(x,t)$ is called as the Lagrange multiplier. A necessary condition for an   optimal solution of model~\eqref{eqn:main_problem} (which is different from our main result of this section) is the stationarity of this Lagrangian functional; i.e., $\displaystyle \lim_{\varepsilon \to 0} \frac{L(u+ \varepsilon \delta u, v, \psi)  - L(u,v,\psi)}{\varepsilon} =0$, where $\delta u$ is an increment of $u$ satisfying equations~\eqref{eqn:difference_2} and~\eqref{eqn:difference_3} for the homogeneous case of $\delta v =0$. This leads us to the following \textit{adjoint problem}. 
\begin{subequations}   \label{eqn:adjoint}
 \begin{eqnarray}
& \psi_{tt}+{(k(x)\psi_{xx})}_{xx}=-2\big[u(x,t;v)-y(x,t)\big], & \forall (x,t)\in \Omega  \label{eqn:adjoint_1} \\
&  \text{with final conditions: } \nonumber \\
& \psi (x,T)=0,\ \ \ \psi_t(x,T)=0,\ \ \ x\in (0,L)\label{eqn:adjoint_2}  \\
&   \text{with boundary conditions:} \nonumber \\
& \psi (0,t)=  \psi (L,t) = \psi_{xx}(0,t) =  \psi_{xx}(L,t)=0 & \forall t \in (0,T). \label{eqn:adjoint_3}
 \end{eqnarray} 
\end{subequations} 

For a given $v\in V_{ad}$, let $\psi (x,t;v)$ denote the weak solution of the problem~\eqref{eqn:adjoint} in the space $\Homega$. By substituting $t$ by $T-t$, we observe that model~\eqref{eqn:adjoint} transforms into a model of the same type as model~\eqref{eqn:euler_bernoulli}. Hence, it follows from  Lemma~\ref{lem:weak_a} that this solution, $\psi$, is also unique. 

Third, we prove an integral identity in the lemma below.

\begin{lemma}\label{lemma2}
 Let  $\delta u(x,t;v)$ and $\psi (x,t;v)$ be the weak solutions for model~\eqref{eqn:difference} and~\eqref{eqn:adjoint}, respectively. Then the following integral identity holds:
$$2\int^T_0{\int^L_0{\bigg[u(x,t;v)-y(x,t)\bigg]\delta u(x,t;v)dx}dt}=-\int^L_0{\psi (x,0;v)\delta v(x)dx}$$
for all $v\in V_{ad}$.
\end{lemma}
\begin{proof}
 Multiplying equation~\eqref{eqn:adjoint_1} by $\delta u$ and integrating on $\Omega$, we have 
\begin{equation}\label{eqn:lemma_integral_1}2 \int^T_0{\int^L_0{\Big[u(x,t;v)-y(x,t)\Big]\delta u(x,t;v)dx}dt}=-\int^T_0{\int^L_0{\Big[\psi_{tt}+{(k(x)\psi_{xx})}_{xx}\Big]\delta udx}dt}.
\end{equation}
\noindent Consider the two terms on the right-hand side of equation~\eqref{eqn:lemma_integral_1}. 
 
 \begin{subequations}\label{eqn:lemma_integral_2}
\begin{eqnarray} 
 \int^T_0{\int^L_0{\psi_{tt} \delta udx}dt} = & \displaystyle  \int^L_0{\bigg[{\psi}_{t} \delta u\Big|_0^T} - \int^T_0{\psi}_{t} (\delta u)_t dt\bigg] dx  \label{eqn:lemma_integral_2_1} &\\
 = &  - \displaystyle \int^L_0{\bigg[{\psi} (\delta u)_t \Big|_0^T} -  \int^T_0{\psi} (\delta u)_{tt} dt\bigg] dx \label{eqn:lemma_integral_2_2} & \\
  = &   \displaystyle \int^L_0{{\psi(x,0;v)} \delta v(x) dx} \label{eqn:lemma_integral_2_3} + \int^T_0{\int^L_0{\psi (\delta u)_{tt} dxdt}}, &
 \end{eqnarray}
\end{subequations}
 and
 \begin{subequations}
 \label{eqn:lemma_integral_3}
\begin{eqnarray}
 &\displaystyle \int^T_0{ \int^L_0 {\big(k(x)\psi_{xx}\big)}_{xx}} \delta u dx dt &=  \displaystyle \int^T_0{\Big[(k(x)\psi_{xx})_x \delta u\Big|_0^L} - \int^L_0{(k(x)\psi_{xx})_x (\delta u)_x dx}\Big]dt \label{eqn:lemma_integral_3_1} \\
 && = \displaystyle - \int^T_0{\Big[k(x)\psi_{xx} (\delta u)_x\Big|_0^L} - \int^L_0{k(x)\psi_{xx} (\delta u)_{xx} dx}\Big]dt \label{eqn:lemma_integral_3_2}  \\
 && = \displaystyle \int^T_0{\Big[k(x)\psi_{x} (\delta u)_{xx} \Big|_0^L} - \int^L_0{\psi_x(k(x)(\delta u)_{xx})_x dx}\Big]dt \label{eqn:lemma_integral_3_3} \\
 && = \displaystyle -  \int^T_0{\Big[\psi(k(x)(\delta u)_{xx})_x \Big|_0^L} - \int^L_0{ \psi\big(k(x)(\delta u)_{xx}\big)_{xx} dx}\Big]dt \label{eqn:lemma_integral_3_4}  \\
  && = \displaystyle \int^T_0{ \psi(k(x)(\delta u)_{xx})_{xx}}dxdt \label{eqn:lemma_integral_3_5}.
\end{eqnarray}
\end{subequations}
In equation~\eqref{eqn:lemma_integral_2_1}, we integrate the left-hand side by parts. Equation~\eqref{eqn:lemma_integral_2_2} follows (i) since the first term on the right-hand side of  equation~\eqref{eqn:lemma_integral_2_1} is zero from the initial condition~\eqref{eqn:difference_2}  and the condition~\eqref{eqn:adjoint_2}, and (ii) by integrating the second term of the right-hand side of equation~\eqref{eqn:lemma_integral_2_1} again by parts.  Equation~\eqref{eqn:lemma_integral_2_3} follows from the initial condition~\eqref{eqn:difference_2}  and the final condition~\eqref{eqn:adjoint_2}.  Similarly, in equation~\eqref{eqn:lemma_integral_3}, we integrate by parts. Then, the first term on the right-hand side of each of the equations~\eqref{eqn:lemma_integral_3_1}-\eqref{eqn:lemma_integral_3_4} is zero from the boundary conditions~\eqref{eqn:difference_3}  and~\eqref{eqn:adjoint_3} while the second term  follows directly from the preceding equation by integration by parts.

Substituting the expressions in equations~\eqref{eqn:lemma_integral_2_3} and~\eqref{eqn:lemma_integral_3_5} in the right-hand side of equation~\eqref{eqn:lemma_integral_1}, we have
\begin{eqnarray} & \label{eqn:lemma_integral_4}
\displaystyle 2 \int^T_0{\int^L_0{\Big[u(x,t;v)-y(x,t)\Big]\delta u(x,t;v)dx}dt} & =  \nonumber \\    
&  & \hspace{-1.5in} - \int^T_0 {\int^L_0 \Big[{(\delta u)_{tt} }  + (k(x)(\delta u)_{xx})_{xx}\Big]\psi}dxdt - \int^L_0{{\psi(x,0;v)} \delta v(x) dx}.
\end{eqnarray}
From equation~\eqref{eqn:difference_1}, the first term on the right-hand side of equation~\eqref{eqn:lemma_integral_4} is zero. This completes the proof.
\end{proof}

With these three results, we derive the Fr\'echet derivative of the cost functional of model~\eqref{eqn:main_problem}.

\begin{lemma}\label{lem:frechet}
 The cost functional given by \eqref{eqn:main_problem_1} is Fr\'echet differentiable and has a derivative given by $J'_{\alpha}(v)= - \psi(x,0)+2\alpha v$,
where $\psi$ is the weak solution of the adjoint problem \eqref{eqn:adjoint}.
\end{lemma}
\begin{proof}
Consider an increment $\delta v \in \Vad$ in the cost functional $J_{\alpha}(v)$ such that $\delta J_{\alpha}(v)=J_{\alpha}(v+\delta v)-J_{\alpha}(v) $. Then, the Fr\'echet derivative $J'_{\alpha}(v)$  is a function $g$ satisfying the following:
\begin{equation} \label{eqn:frechet_def}
 \delta J_{\alpha}(v)={\big \langle g,\delta v \big \rangle}_{\Ltwo} + o\bigg(\norm{\delta v}_{\Ltwo}^2\bigg).
\end{equation}
We have:
\begin{subequations}\label{lem:frechet_1}
\begin{align}
\delta J_{\alpha}(v)&=J_{\alpha}(v+\delta v)-J_{\alpha}(v) \label{lem:frechet_1_1}\\
&=\int_0^T\int_0^L 2\Big[u(x,t;v)-y(x,t)\Big]\delta u(x,t)dxdt+\int_0^T\int_0^L \Big[\delta u(x,t)\Big]^2dxdt \nonumber\\
& \, \, + 2\alpha {\big \langle v,\delta v \big \rangle}_{\Ltwo}+\alpha \norm{\delta v}_{\Ltwo}^2. \label{lem:frechet_1_2} \\
& = -\int _0^L\psi(x,0)\delta v(x)dx + 2\alpha {\langle v,\delta v\rangle}_{\Ltwo}+o\bigg(\norm{\delta v}_{\Ltwo}^2\bigg). \label{lem:frechet_1_3}\\
& =  \langle -\psi(x,0) + 2 \alpha v,\delta v\rangle_{\Ltwo}+o\bigg(\norm{\delta v}_{\Ltwo}^2\bigg). \label{lem:frechet_1_4}
\end {align}
\end{subequations}

Here  equation~\eqref{lem:frechet_1_2} follows from the definition of $J_\alpha$ in equation~\eqref{eqn:main_problem_1}, the definition of $\delta u$,  Definition~\ref{def:ltwo}, and by using the identity $a^2-b^2= (a+b)(a-b)$. Then, equation~\eqref{lem:frechet_1_3} follows by: (i) substituting the first term on the right-hand side of equation~\eqref{lem:frechet_1_2} by that from Lemma~\ref{lemma2}, and (ii)  by bounding the second term on the right-hand side of equation~\eqref{lem:frechet_1_2} by that from  Proposition~\ref{cor:frechet}. Equation~\eqref{lem:frechet_1_4} follows from Definition~\ref{def:ltwo}. The result follows from equation~\eqref{eqn:frechet_def}.
\end{proof}

\begin{theorem} \label{thm:necessarycondition}
Let $u^*$ be the weak solution of equations~\eqref{eqn:main_problem_2}-\eqref{eqn:main_problem_4} corresponding to a $v^*$. Further, let $\psi^*$ be the weak solution of the adjoint problem given by ~\eqref{eqn:adjoint} corresponding to  $u^*$. Then, a necessary condition for $v^*$ to be the optimal solution for model~\eqref{eqn:main_problem} is  
\begin{equation}\label{eqn:lnecessarycondition}
   \big \langle J'_{\alpha}(v^*),v-v^* \big \rangle_{\Ltwo} = \big \langle - \psi^*(x,0)+2\alpha v^*,v-v^* \big \rangle_{\Ltwo}  \ge 0, \quad \forall  v \in \Vad.
\end{equation}
\end{theorem}

\begin{proof}
Using the formula for $ J'_{\alpha}(v^*)$ from Lemma~\ref{lem:frechet}, the result follows from~\cite{Zeidler1985}.
\end{proof}

\section{A gradient-type method for a numerical solution of model~\eqref{eqn:main_problem}} \label{sec:lipschitz}
In this section, we provide an iterative method for a numerical solution of model~\eqref{eqn:main_problem}; further, we show the existence of a Lipschitz constant for the Fr\'echet derivative of the cost functional. We recall that the Fr\'echet  differentiability of the cost functional given by~\eqref{eqn:main_problem_1} is known from Lemma~\ref{lem:frechet}. 

In Algorithm~\ref{alg:algorithm} we provide a pseudocode for a gradient-based  descent method  that seeks to minimize the cost functional $J_{\alpha}(v)$. We start with an initial control, $v^0$, from the set of admissible controls and iteratively update it to approximate an optimal solution within a chosen tolerance $\varepsilon$. The weak solutions in Step~\ref{alg:u} and Step~\ref{alg:psi} can be computed using the Galerkin method, which is a widely used approach for obtaining weak solutions of hyperbolic partial differential equations.  This method rests on constructing an appropriate sequence of approximations that ultimately converge to the weak solution; see, e.g.,\cite{Baysal2019, Ladyzhenskaya2013} for details. Step~\ref{alg:j} is accomplished via Lemma~\ref{lem:frechet}. In  Step~\ref{alg:update}, the parameter $\beta$ is a step-size that may be chosen as is standard in any descent-based optimization method, see, e.g.,~\cite{Ritter2023}. A common choice for the step-size that ensures convergence is the inverse of the Lipschitz constant, see, e.g~\cite[page~14]{Hasan2009}. This choice ensures the update condition in Step~\ref{alg:update} is satisfied. Further, in  the updating rule of Step~\ref{alg:update},  the term $J'_{\alpha}(v^k)$ plays the
role of the steepest descent direction for the minimization. If $J'_{\alpha}(v^k)$ is zero at any iteration $k$, we terminate the algorithm since $v^k$ is a stationary point of model~\eqref{eqn:main_problem}; however,  this does not guarantee that $v^k$ is an optimal solution. 

 If the gradient of the functional $J_{\alpha}$ is Lipschitz continuous, then the parameter $\beta_k$ can be selected from the following condition (see, e.g.,~\cite{Polyak1987}): 
 \begin{equation} \label{eqn:lipschitz}
 \epsilon_0\le \beta_k \le \frac{2}{\mathcal{L}+2\epsilon_1},
\end{equation}
where $\mathcal{L}>0$ is the Lipschitz constant, and  $\epsilon_0, \epsilon_1>0$ are constant parameters. In particular, setting $\epsilon_0=\frac{1}{\mathcal{L}}$ and $ \epsilon_1=\frac{\mathcal{L}}{2}$ we obtain a constant step-size parameter $\beta_k=\frac{1}{\mathcal{L}}$. The following lemma derives such a Lipschitz constant. 

\begin{algorithm}[!htb]
\caption{Sketch of an algorithm for a numerical solution of model~\eqref{eqn:main_problem}.}\label{alg:algorithm}
\begin{algorithmic}[1]
\Require An initial solution $v^0(x) \in \Vad$; an appropriate step-size parameter $\beta^k$ for iteration $k$; a tolerance $\varepsilon$. 
\Ensure  An optimal solution $v^*(x)$ for model~\eqref{eqn:main_problem} subject to the chosen tolerance $\varepsilon$.
\Repeat \, at iteration $k$:
\State \label{alg:u} Given $v^k(x)$, obtain a weak solution, $u^k(x,t)$, for  model~\eqref{eqn:euler_bernoulli}.
\State \label{alg:psi} Given this $u^k(x,t)$, obtain a weak solution $\psi^k\big(x,t;v^k(x)\big)$ for model~\eqref{eqn:adjoint}. 
\State \label{alg:j} Given this $\psi^k\big(x,t;v^k(x)\big)$, obtain the gradient $J'_{\alpha}\big(v^k(x)\big)$.
\State \label{alg:update} Choose $\beta_k$ such that $J_{\alpha}\big(v^{k+1}(x)\big)<J_{\alpha}\big(v^k(x)\big)$. Update $v$: $v^{k+1}(x) \gets v^k(x) - \beta^k J'_{\alpha}(v^k(x))$.  
 \Until{$\norm{v^{k+1}(x)-v^k(x)}_{\Ltwo} \le \varepsilon$}; $v^{k+1}(x) \gets v^*(x)$. 
\end{algorithmic}
\end{algorithm}

\begin{lemma} \label{lem:lipschitz} 
Let the function $k(x)$ be positive and finite-valued with $k(x) \in L^\infty(0,L)$, $w(x) \in H^2(0,L)$ and $F(x,t) \in \Lomega$. 
\begin{enumerate}[label={L.(\roman*)}] 
\item \label{lem:lipschitz_a} The following estimate holds: 
\[ \norm{J'_{\alpha}(v +\delta v)-J'_{\alpha}(v)} _{\Ltwo}^2 =\int_0^L \big(-\delta \psi(x,0;v)+2\alpha \delta v\big)^2dx. \]
\item \label{lem:lipschitz_b} There exists a Lipschitz constant, $\mathcal{L}$, such that
\[ \label{eqn:lipschitz}
       \norm{J'_{\alpha}(v+\delta v)-J'_{\alpha}(v)} _{\Ltwo} \le \mathcal{L} \norm{\delta v}_{\Ltwo}. \]
\end{enumerate}
\end{lemma}

\begin{proof} 
The function $\delta \psi(x,t;v)=\psi(x,t;v+\delta v)-\psi(x,t;v)$ is the solution of the following backward hyperbolic problem:
\begin{subequations}   \label{eqn:adjointdifference}
 \begin{align}
& \delta \psi_{tt}+{\big(k(x)\delta \psi_{xx}\big)}_{xx}=-2\big[\delta u(x,t;v)\big], & \forall (x,t)\in \Omega  \label{eqn:adjointdif_1} \\
& \, \,  \text{with final conditions:}  \nonumber\\
& \delta \psi (x,T)=0, \, \delta \psi_t(x,T)=0, & \forall x\in (0,L)\label{eqn:adjointdif_2}  \\
& \, \,  \text{with boundary conditions:} \nonumber \\
& \delta \psi (0,t)= \delta \psi (L,t) = \delta \psi_{xx}(0,t) = \delta \psi_{xx}(L,t)=0 & \forall t \in (0,T). \label{eqn:adjointdif_3}
 \end{align} 
\end{subequations}

Replacing $t$ by $T-t$ in model~\eqref{eqn:adjointdifference}, we note that model~\eqref{eqn:adjointdifference} with  homogeneous initial data is of the same type as model~\eqref{eqn:euler_bernoulli}. Hence,
it follows from Lemma~\ref{lem:weak} that there exists a constant $C_0$ such that
\begin{equation} \label{eqn:lipschitz_co}
\norm{\delta \psi}_{H^{2,1}(\Omega)} \le 2C_0 \norm{\delta u}_{\Ltwo}.
\end{equation}
We have from Definition~\ref{def:htwo} and equation~\eqref{eqn:lipschitz_co} that
\begin{equation}
 \label{eqn:lipschitz_co_2_1}
 \norm{\delta \psi}_{\Lomega}^2+\norm{\delta \psi_t}_{\Lomega}^2 \le 4C_0^2 \norm{\delta u}_{\Ltwo}^2,
 \end{equation}
while using Proposition~\ref{cor:frechet} it further follows that 
\begin{equation}
\label{eqn:lipschitz_co_2_1}
\norm{\delta \psi}_{\Lomega}^2+\norm{\delta \psi_t}_{\Lomega}^2 \le 4C_0^2 C_1^2 \norm{\delta v}_{\Ltwo}^2.  
 \end{equation}
Consider the following identity.
\begin{equation} \label{eqn:deltapsi_0}
 \big[\delta \psi(x,0)\big]^2 = \Big (-\int_0^t  \delta \psi_t(x,\tau)d\tau+\delta \psi(x,t)\Big)^2, \, \, \forall t \in (0,T).
\end{equation}
We then have:
\begin{subequations} \label{eqn:deltapsi_0_estimate}
 \begin{eqnarray}
 & \bigg[\delta \psi(x,0)\bigg]^2  & \le 2\Big ( \int_0^t  \delta \psi_t(x,\tau)d\tau\Big )^2+2\big(\delta \psi(x,t)\big)^2.  \label{eqn:deltapsi_0_estimate_1} \\
\implies & \displaystyle \bigg[\delta \psi(x,0)\bigg]^2 & \le  2t\int_0^t  (\delta \psi_t(x,\tau))^2d\tau+2(\delta \psi(x,t))^2.   \label{eqn:deltapsi_0_estimate_2} \\
\implies & \displaystyle \bigg[\delta \psi(x,0)\bigg]^2  &  \le   T\int_0^T  (\delta \psi_t(x,t))^2dt+\frac{2}{T}\int^T_0(\delta \psi(x,t))^2dt \label{eqn:deltapsi_0_estimate_3} \\
 \implies & \displaystyle \int^L_0[\delta \psi(x,0)]^2dx  &  \le   T\int^L_0\int_0^T  (\delta \psi_t(x,t))^2dtdx+\frac{2}{T}\int^L_0\int_0^T (\delta \psi(x,t))^2dtdx .\label{eqn:deltapsi_0_estimate_4} \\
 \implies & \displaystyle \int^L_0[\delta \psi(x,0)]^2dx  &  \le   max\Big\{T,\frac{2}{T}\Big\} \Big (\norm{\delta \psi}_{\Lomega}^2+\norm{\delta \psi_t}_{\Lomega}^2 \Big ).\label{eqn:deltapsi_0_estimate_5} \\
  \implies & \displaystyle \int^L_0[\delta \psi(x,0)]^2dx  &  \le   max\Big \{T,\frac{2}{T}\Big \}4C_0^2 C_1^2 \norm{\delta v}_{\Ltwo}^2.\label{eqn:deltapsi_0_estimate_6}
 \end{eqnarray}
\end{subequations}

Here, equation~\eqref{eqn:deltapsi_0_estimate_1} follows from equation~\eqref{eqn:deltapsi_0} and by applying the inequality $(a+b)^2 \le 2(a^2+b^2)$.  
Then, equation~\eqref{eqn:deltapsi_0_estimate_2} follows from the Cauchy-Schwartz inequality applied to the first term on the  right-hand side of equation~\eqref{eqn:deltapsi_0_estimate_1}. Since equation~\eqref{eqn:deltapsi_0_estimate_2} holds for all $t \in (0,T)$ and the second term of equation~\eqref{eqn:deltapsi_0_estimate_2} is always non-negative, equation~\eqref{eqn:deltapsi_0_estimate_3} follows. Equation~\eqref{eqn:deltapsi_0_estimate_4} follows by further integrating both sides of equation~\eqref{eqn:deltapsi_0_estimate_3} on $[0,L]$.   Equation~\eqref{eqn:deltapsi_0_estimate_5} holds from the definition of the norm in space $\Lomega$, while  equation~\eqref{eqn:deltapsi_0_estimate_6} follows from  equation~\eqref{eqn:deltapsi_0}.

Next, from Lemma~\ref{lem:frechet}, we have
\begin{subequations}
 \label{eqn:lemma_difference_of_cost}
\begin{eqnarray}
 &\displaystyle J'_{\alpha}(v +\delta v)-J'_{\alpha}(v) &=  \displaystyle -\psi(x,0;v+\delta v)+2\alpha (v+\delta v)+\psi(x,0;v)-2\alpha v  \\
 && = \displaystyle -\delta \psi(x,0;v)+2\alpha \delta v,  
 \end{eqnarray}
\end{subequations}
where $\delta \psi(x,t;v)=\psi(x,t;v+\delta v)-\psi(x,t;v)$ is the solution of the boundary value problem~\eqref{eqn:adjointdifference}. From Definition~\ref{def:ltwo}, we then have 
 \begin{equation} \label{eqn:norm_of_Frechet}
\norm{J'_{\alpha}(v +\delta v)-J'_{\alpha}(v)} _{\Ltwo}^2 =\int_0^L \big(-\delta \psi(x,0;v)+2\alpha \delta v\big)^2dx,
\end{equation}
which proves Lemma~\ref{lem:lipschitz_a}. Applying the inequality $(a+b)^2 \le 2(a^2+b^2)$ to equation~\eqref{eqn:norm_of_Frechet}, we have 
\begin{equation}  \label{eqn:norm_of_Frechet_2}
\norm{J'_{\alpha}(v +\delta v)-J'_{\alpha}(v)} _{\Ltwo}^2 \le 2 \int_0^L \big(\delta \psi(x,0;v))^2dx + 8\alpha^2 \int_0^L \delta v^2dx.
\end{equation}

From equation~\eqref{eqn:deltapsi_0_estimate_6} and equation~\eqref{eqn:norm_of_Frechet_2} we obtain an estimate of the Lipschitz constant of Lemma~\ref{lem:lipschitz_b} as $\mathcal{L}=\Big (8 \max \big \{T,\frac{2}{T}\big \} C_0^2 C_1^2 +8\alpha ^2\Big)^{\frac{1}{2}}$.
\end{proof}

\section{Conclusions} \label{sec:conclusion}
In this work, we formulate an optimal control problem via the initial velocity of the equation of hyperbolic type within the framework of Euler–Bernoulli beam theory. By obtaining energy estimates for the corresponding initial-boundary value problem, we establish the existence and uniqueness of the weak solution, $v(x)$, of the control function.  We construct an adjoint problem corresponding to the optimal control formulation and derive an explicit expression for the Fr\'echet gradient of the cost functional. Further, we prove the Lipschitz continuity of this gradient which allows us to estimate the step-size parameter in descent-based optimization algorithms using the associated Lipschitz constant. These results enable an efficient numerical computation of optimal solutions via a minimizing sequence. Future work may examine control problems where the initial velocity belongs to function spaces beyond $\Ltwo$, e.g., $H^2(0,L)$.

\section*{Funding}
We thank T{\"U}B\.{I}TAK-BIDEB 2221 – Fellowships for Visiting Scientists and Scientists on Sabbatical Leave Support Program (2024/1) that partially supported us. We further acknowledge the support of the Heilbronn Institute for Mathematical Research in enabling this work. 


\bibliographystyle{chicago}
\bibliography{thebib2}

\end{document}